\catcode`\@=11
\magnification 1200
\hsize=140mm \vsize=200mm
\hoffset=-4mm \voffset=-1mm
\pretolerance=500 \tolerance=1000 \brokenpenalty=5000

\catcode`\;=\active
\def;{\relax\ifhmode\ifdim\lastskip>\z@
\unskip\fi\kern.2em\fi\string;}

\catcode`\:=\active
\def:{\relax\ifhmode\ifdim\lastskip>\z@\unskip\fi
\penalty\@M\ \fi\string:}

\catcode`\!=\active
\def!{\relax\ifhmode\ifdim\lastskip>\z@
\unskip\fi\kern.2em\fi\string!}\catcode`\?=\active
\def?{\relax\ifhmode\ifdim\lastskip>\z@
\unskip\fi\kern.2em\fi\string?}

\frenchspacing

\newif\ifpagetitre        \pagetitretrue
\newtoks\hautpagetitre    \hautpagetitre={\hfil}
\newtoks\baspagetitre     \baspagetitre={\hfil}
\newtoks\auteurcourant    \auteurcourant={\hfil}
\newtoks\titrecourant     \titrecourant={\hfil}
\newtoks\hautpagegauche   \newtoks\hautpagedroite
\hautpagegauche={\hfil\tensl\the\auteurcourant\hfil}
\hautpagedroite={\hfil\tensl\the\titrecourant\hfil}

\newtoks\baspagegauche
\baspagegauche={\hfil\tenrm\folio\hfil}
\newtoks\baspagedroite
\baspagedroite={\hfil\tenrm\folio\hfil}

\headline={\ifpagetitre\the\hautpagetitre
\else\ifodd\pageno\the\hautpagedroite
\else\the\hautpagegauche\fi\fi}

\footline={\ifpagetitre\the\baspagetitre
\global\pagetitrefalse
\else\ifodd\pageno\the\baspagedroite
\else\the\baspagegauche\fi\fi}

\font\twbf=cmbx12\font\sc=cmcsc10

\def\date{le\ {\the\day}\
\ifcase\month\or Janvier\or \F\'evrier\or Mars\or Avril
\or Mai\or Juin\or Juillet\or Ao\^ut\or Septembre
\or Octobre\or Novembre\or D\'ecembre\fi\ {\the\year}}

\def\cf{{\it cf.\/}\ }    \def\ie{{\it i.e.\/}\ }
 \def\up#1{\raise 1ex\hbox{\sevenrm#1}}
\def\cqfd{\unskip\kern 6pt\penalty 500
\raise -2pt\hbox{\vrule\vbox to 10pt{\hrule width 4pt\vfill
\hrule}\vrule}\par}
\catcode`\@=12
\def \bg {\bigskip \goodbreak}
\def \sn {\nobreak \smallskip}

\def\ref#1&#2&#3&#4&#5\par{\par{\leftskip = 5em {\noindent
\kern-5em\vbox{\hrule height0pt depth0pt width
5em\hbox{\bf[\kern2pt#1\unskip\kern2pt]\enspace}}\kern0pt}
{\sc\ignorespaces#2\unskip},\
{\rm\ignorespaces#3\unskip}\
{\sl\ignorespaces#4\unskip\/}\
{\rm\ignorespaces#5\unskip}\par}}

\def\exo#1{\goodbreak\vskip 12pt plus 20pt minus 2pt
\line{\noindent\hss\bf
\uppercase\expandafter{\romannumeral#1}\hss}\nobreak\vskip
12pt }
\def \titre#1\par{\null\vskip
1cm\line{\hss\vbox{\twbf\halign
{\hfil##\hfil\crcr#1\crcr}}\hss}\vskip 1cm}

\def \frac#1#2{{{#1}\over {#2} }}

\def \reel{ {\rm I}\!{\rm R}}

\def \rat{ {\rm Q}\kern-.65em {}^{{}_/ }}

\def \adh {\overline}

\def\N#1{\left\Vert #1\right\Vert }

\def\dess#1by#2(#3){\vbox to #2{
\hrule width #1 height 0pt depth 0pt
\vfill\special{picture #3}
}}

\baselineskip=18pt
\hfuzz=0,5cm
\def\adh{\overline}
\centerline {Sur les op\'erateurs factorisables par $OH$}
\centerline {par Gilles Pisier}
\def\om{\otimes_{\min}}
Note pr\'esent\'ee par

\bg
\bg
\bg
\underbar {R\'esum\'e}. Soient $H,K$ deux espaces de Hilbert. Soient $E
\subset
B(H)$ et $F \subset B(K)$ deux espaces de Banach d'op\'erateurs, au
sens de
[1,2]. On \'etudie les op\'erateurs $u : E \to F$ qui admettent une
factorisation
$E \to OH \to F$ avec des applications compl\`etement born\'ees \`a
travers l'espace
de Hilbert d'op\'erateurs $OH$ que nous avons introduit et \'etudi\'e
dans une note
pr\'ec\'edente. Nous donnons une caract\'erisation de ces op\'erateurs
qui permet de
faire une th\'eorie enti\`erement analogue au cas des op\'erateurs
entre espaces de
Banach qui se factorisent par un Hilbert.
\bg
\bg
\underbar {English Abstract}. Let $H,K$ be Hilbert spaces. Let $E
\subset B(H)$
and $F \subset B(K)$ be operator spaces in the sense of [1,2]. We study
the
operators $u : E \to F$ which admit a factorization $E \to OH \to F$
with
completely bounded maps through the operator Hilbert space $OH$ which
we have
introduced and studied in a recent note. We give a characterization of
these
operators which allows to develop a theory entirely analogous to that
of
operators between Banach spaces which can be factored through a Hilbert
space.
\bg
\bg
\bg
\underbar {Abridged English version}. In this note we continue the
study of the
operator Hilbert space $OH$ introduced in our previous note [6]. Let
$E,F$ be
two operator spaces in the sense of [1,2].
We can assume $E\subset B(H)$ and $F\subset B(K)$ for some Hilbert
spaces $H$
and $K$. We will denote by $H \otimes
K$ the Hilbertian tensor product. We denote by $E \otimes_{\min} F$ the
minimal (or spatial) tensor product in $B(H \otimes K)$.    We consider
the
space $\Gamma_{oh}(E,F)$ of all operators $u : E \to F$ for which there
are an
index set $I$ and completely bounded (in short c.b.) maps $A : OH(I)
\to F$ and
$B : E \to OH(I)$ such that $u = AB$. In that case, we say that $u$
factors
through $OH$. We denote by $\gamma_{oh}(u)$ the infimum of $\N {A}_{cb}
\N
{B}_{cb}$ over all possible such factorizations. If $v\in E \otimes F$
we denote by $\gamma_{oh}(v)$ the above norm for the associated
operator
from $E^*$ to $F$ (or from $F^*$ to $E$).
Note that
$u$ factors through $OH$ iff its adjoint $u^*$ also does and
$\gamma_{oh}(u)=\gamma_{oh}(u^*)$.
 If $F$ is the antidual $\adh {E^*}$
(\ie the dual with the conjugate complex multiplication) equipped with
its
operator space structure in the sense of [1,2] we say that a map $u : E
\to \adh
{E^*}$ is positive if the associated sesquilinear form is positive \ie
if $u(x)
(x) \geq 0$ for all $x$ in $E$.

\proclaim THEOREM 1. If $u : E \to \adh {E^*}$ is positive and
completely
bounded, then $u \in \Gamma_{oh} (E,\adh {E^*})$ and $\gamma_{oh} (u) =
\N
{u}_{cb}$. Moreover, every map $u$ in $\Gamma_{oh} (E,\adh {E^*})$ can
be
written as $u_1 - u_2 + i(u_3 - u_4)$ with $u_1,...,u_4$ positive and
such that
$\gamma_{oh} (u_j) \leq \gamma_{oh}(u)$ for all $j = 1,...,4$.

\underbar {Remark}. Contrary to the Banach space case where
Grothendieck's
theorem says that every bounded map $u : C \to \adh {C^*}$ factors
through
$\ell_2$, it is not true that every c.b. map $u : B(H) \to \adh
{B(H)^*}$
factors through $OH$.

\underbar {Definition}. Let $E$ be an operator space and let $F$ be a
Banach
space. An operator $u : E \to F$ will be called $(2,oh)$-summing if
there is a
constant $C$ such that for all finite sequences $(x_i)$ in $E$ we have
$$(\sum \N {u(x_i)}^2)^{1/2} \leq C \N {\sum x_i \otimes \adh
{x_i}}^{1/2}_{E\otimes_{\min} \adh E}.$$ We denote by $\pi_{2,oh}(u)$
the
smallest constant $C$ for which this holds. (Similar but different
spaces
have already been considered in [1,3].) There is of course an analogue
of the
``Pietsch factorization" (\cf e.g. [7]) for these operators, but the
main point
is that they provide a convenient description of the dual tensor norm
to the norm
$\gamma_{oh}$, in ``complete" analogy with the Banach space case, as
follows.

\proclaim THEOREM 2. Let $E,F$ be operator spaces and let $C > 0$ be a
constant. The following properties of a map $u : E \to F$ are
equivalent.\sn
(i) $u \in \Gamma_{oh}(E,F)$ and $\gamma_{oh}(u) \leq C$.\sn
(ii) For all $v : F \to \ell_2$ such that $\pi_{2,oh} (v) \leq 1$
the composition $vu$ can be factorized as $vu = AB$ with
$B : E \to OH$ completely bounded and
$A : OH \to \ell_2$ Hilbert Schmidt satisfying $\N {A}_{HS} \N
{B}_{cb}\leq C$.

 From the preceding it is easy to deduce a
  description of the   tensor norm   $\gamma_{oh}^*$ which is dual to
  the norm
$\gamma_{oh}$, exactly as in the Banach space case, as follows.

\proclaim COROLLARY 3. Let $E_1,E_2$ be two   operator spaces. Assume
(to
simplify the statement) $E_2$ reflexive. For any  operator $u$ from
$E_1$ to
$E_2$ let $$\gamma_{oh}^*(u)=\sup\{ |<u,v>|\quad  |\quad v\in
E_1\otimes
E_2^*,\  \gamma_{oh}(v)\leq 1\}.$$ We have then
$$\gamma_{oh}^*(u)=\inf\{ \pi_{2,oh}(B) \pi_{2,oh}(A^*) \}$$
where the infimum runs over all possible  factorizations of $u$ of the
form $u=AB$ with operators $B: E_1\to \ell_2$  and $A:\ell_2\to E_2$
such
that $B$ and $A^*$ are $(2,oh)-$summing.

Actually, a more general class of tensor norms which we called
``$\gamma$-norms" in [7] can be treated in very much the same way as
above for
$\gamma _{oh}$. These results allow the development of a theory of type
and
cotype or of a ``local theory" (see e.g. [7] for all this) in the
category of
operator spaces.

I am very grateful to David Blecher and Vern Paulsen for stimulating
conversations on the subject of this note and the preceding one [6].
\bg
\bg
\bg
Cette note fait suite \`a la note pr\'ec\'edente [6] et annonce les
r\'esultats d'un
article \`a paratre.

Nous renvoyons \`a [1,2] pour la th\'eorie des espaces d'op\'erateurs
et \`a notre
travail [6] pour tout ce qui concerne l'espace $OH$ (resp. $OH(I))$ qui
est
l'analogue de $\ell_2$ (resp. $\ell_2(I))$ dans la cat\'egorie des
espaces
d'op\'erateurs.

Soient $H,K$ deux Hilbert, $E \subset B(H),F \subset B(K)$ deux
sous-espaces
ferm\'es.
On notera $E \otimes_{\min}F$ le produit tensoriel compl\'et\'e pour la
norme
induite par l'espace $B(H \otimes K)$ des op\'erateurs born\'es sur le
produit
tensoriel hilbertien $H \otimes K$. Nous renvoyons \`a
[1] ou [2] pour la d\'efinition du dual $E^*$ d'un espace d'op\'erateur
$E$
ainsi que pour la notion d'application compl\`etement born\'ee.

Nous dirons qu'un op\'erateur $u : E \to F$ se factorise par $OH$ s'il
existe un
ensemble $I$ et des applications compl\`etement born\'ees $B : E \to
OH(I)$ et $A :
OH(I) \to F$ telles que $u = AB$. On pose $\gamma_{oh}(u) = \inf \{\N
{A}_{cb}
\N {B}_{cb}\}$ o\`u l'infimum porte sur toutes les factorisations
possibles de
$u$. On notera $\Gamma_{oh} (E,F)$ l'espace des applications $u : E \to
F$ qui
se factorisent par $OH$. C'est un espace de Banach quand on le munit de
la
norme $\gamma_{oh}$. Nous noterons $E^*$ (resp. $\adh {E^*})$ le dual
(resp.
l'antidual) de $E$ au sens des espaces d'op\'erateurs (\cf [1,2]), de
sorte
qu'une application lin\'eaire $u : E \to \adh {E^*}$ correspond \`a une
application
sesquilin\'eaire sur $E \times  E$. Nous dirons que $u : E \to \adh
{E^*}$ est
positive si la forme sesquilin\'eaire associ\'ee est positive, \ie si
$u(x)(x) \geq
0~~~~\forall~x \in E$.

\proclaim THEOREME 1. Soit $E \subset B(H)$. Soit $u : E \to \adh
{E^*}$.\sn
(i) Si $u$ est positive alors $u \in \gamma_{oh}(E,\adh {E^*})$ et
$\gamma_{oh}(u) = \N {u}_{cb}$.\sn
(ii) Toute application $u \in \Gamma_{oh} (E,\adh {E^*})$ peut
s'\'ecrire $u_1 -
u_2 + i(u_3 - u_4)$ avec $u_1, u_2, u_3, u_4$ positives de $E$ dans
$\adh
{E^*}$ et telles que
$$\forall~j = 1,...,4~~~~~~~\N {u}_{cb} = \gamma_{oh}(u_j) \leq
\gamma_{oh}(u).$$

Ce th\'eor\`eme signifie que l'espace $\Gamma_{oh}(E,\adh {E^*})$
coincide avec
l'ensemble des combinaisons lin\'eaires d'applications positives
compl\`etement
born\'ees de $E$ dans $\adh {E^*}$.

Parmi les propri\'et\'es \'el\'ementaires des
op\'erateurs appartenant \`a $\Gamma_{oh} (E,F)$, citons les suivantes
: toute
ultraproduit d'applications $u_i$ avec $u_i \in \Gamma_{oh} (E_i,F_i)$
et
$$\sup_{i \in I} \gamma_{oh} (u_i) < \infty$$
se factorise par $OH$ et l'op\'erateur $u$ r\'esultant de
l'ultraproduit v\'erifie
$\gamma_{oh} (u) \leq \sup_{i \in I} \gamma_{oh} (u_i)$. La borne
inf\'erieure
dans la d\'efinition de $\gamma_{oh}(u)$ est atteinte. De plus, si $E_1
\subset
E$ et $F_1 \subset F$ sont des sous-espaces ferm\'es et si on note $q :
E \to
E/E_1$ et $j : F_1 \to F$ les morphismes canoniques alors on a pour
tout
$u : E/E_1 \to F_1$
 $$\gamma_{oh} (u) = \gamma_{oh}(j u q).$$

\underbar {D\'efinition}. Soit $E$ un espace d'op\'erateurs et $F$ un
espace de
Banach. Nous dirons qu'un op\'erateur $u : E \to F$ entre espaces
d'op\'erateurs
est $(2,oh)$-sommant s'il existe une constante $C$ telle que
$$\forall~n  ~~~\forall~x_i \in E~~~~~~(\sum \N {u(x_i)}^2)^{1/2} \leq
C \N {\sum x_i \otimes \adh
{x_i}}^{1/2}_{E\otimes_{\min} \adh E}.$$
On rappelle que l'on a (voir [6])
$$\N {\sum x_i \otimes \adh
{x_i}}^{1/2}_{E\otimes_{\min} \adh E}=\N {\sum x_i \otimes
T_i} _{E\otimes_{\min} OH}$$
o\`u $(T_i)$ est une base orthonormale fix\'ee quelconque de
l'espace $OH$. On notera $\pi_{2,oh}(u)$ la plus petite
constante $C$ v\'erifiant cette propri\'et\'e et $\Pi_{2,oh}
(E,F)$ l'espaces des op\'erateurs $(2,oh)$-sommants de $E$
dans $F$. C'est un espace de Banach muni de la norme
$\pi_{2,oh}$. Cet espace est stable par composition \`a
droite (resp. \`a gauche) par des applications compl\`etement
born\'ees (resp. born\'ees). D'autres espaces du m\^eme genre
(mais diff\'erents) ont d\'ej\`a \'et\'e consid\'er\'es dans [1] et
[3].

\proclaim THEOREME 2. Soient $E,F$ deux espaces d'op\'erateurs. Soit
$C$ une
constante positive. Les propri\'et\'es suivantes d'un op\'erateur $u :
E \to F$ sont
\'equivalentes :\sn
(i) $u \in \Gamma_{oh} (E,F)$ et $\gamma_{oh} (u) \leq C$.\sn
(ii) Pour tout op\'erateur $v \in \Pi_{2,oh} (F, \ell_2)$ l'op\'erateur
$(vu)^* :
\ell_2 \to E^*$ admet une factorisation par $OH$ de la forme $(vu)^* =
wV$ avec
$V : \ell_2 \to OH$ de Hilbert Schmidt et $w : OH \to E^*$
compl\`etement born\'e
tels que
$$\N {V}_{HS} \N {w}_{cb} \leq C \pi_{2,oh} (v).$$

Le th\'eor\`eme pr\'ec\'edent permet de donner (voir  ci-dessous) une
description de la norme duale de la norme $\gamma_{oh}$ enti\`erement
analogue au
cas des normes $\gamma_2$ et $\gamma_2^*$ dans le cadre des espaces de
Banach
(voir par exemple [7]).

Soient $E_1,E_2$ deux espaces d'op\'erateurs. Pour tout $v\in
(E_1\otimes_{\min} E_2
)^*$ on note $I(v)$ la norme de $v$ dans le dual de $E_1\otimes_{\min}
E_2$.
C'est l'analogue de la norme int\'egrale pour les espaces
d'op\'erateurs. La
th\'eorie des op\'erateurs "int\'egraux"   (et des op\'erateurs
"nucl\'eaires") dans ce
nouveau cadre est faite dans [4].

\proclaim PROPOSITION 3.  Soit $E$ un espace d'op\'erateurs. Un
op\'erateur $u:E\to
\ell_2$ est $(2,oh)-$sommant avec $\pi_{2,oh} (u)\le 1$ si et seulement
si il
existe $v\in (E\otimes_{\min} \adh E)^*$ avec $I(v)\le 1$ tel que pour
tout $x$
dans $E$ on a  $v(x\otimes \adh x)\in \reel$ et
$$  \| u(x)\|^2 \leq v(x\otimes \adh x).$$

 En fait on peut placer ces r\'esultats dans
un cadre beaucoup plus large, celui des $\gamma$-normes d\'ej\`a
\'etudi\'es dans [8].
Dans le reste de cette note nous esquissons cette th\'eorie (que nous
d\'evelopperons dans une publication ult\'erieure) dans le style des
id\'ees
originales de [5].

Soit $E$ un espace de Banach.

On notera $B_+(E)$ l'ensemble des \'el\'ements positifs de $E \otimes
\adh {E}$,
c'est-\`a-dire l'ensemble des \'el\'ements qui d\'efinissent une forme
sesquilin\'eaire
de rang fini $\sigma(E^*,E)$-continue et positive sur $E^*\times E^*$.
On note
que $ B_+(E) - B_+(E)$ peut \^etre identifi\'e au sous-espace de $E
\otimes \adh {E}$
form\'e des tenseurs sym\'etriques.

Soit $\gamma : B_+(E) \to \reel_+$ une application additive
positivement
homog\`ene et telle que $0 \leq u \leq v$ implique $0 \leq \gamma (u)
\leq \gamma
(v)$. Nous appellerons ``poids" une telle application $\gamma$. Nous
dirons que
le poids $\gamma$ est raisonnable s'il existe deux constantes $c > 0$
et $C$
telles que pour tout $x$ dans $E$ on a
$$c \N {x}^2 \leq \gamma (x \otimes \adh {x}) \leq C \N {x}^2.$$
Tout poids homog\`ene $\gamma$ donn\'e sur $B_+(E) \subset E \otimes
\adh {E}$ peut
\^etre prolong\'e sur $E \otimes \adh {E}$ en posant
$$\forall~u \in E \otimes \adh {E}~~~~\gamma(u) = \inf \{\gamma(\sum
x_i
\otimes \adh {x_i})^{1/2} \gamma(\sum y_i \otimes \adh {y_i})^{1/2}\}$$
o\`u l'infimum porte sur toutes les repr\'esentations $$u = \sum_1^n
x_i \otimes
\adh {y_i},\quad
x_i \in E,\quad y_i \in E.$$
Si $\gamma$ est raisonnable, ce prolongement d\'efinit une norme sur $E
\otimes
\adh {E}$. Plus g\'en\'eralement, soit $E_1,E_2$ deux espaces de Banach
et soit
$\gamma_1,\gamma_2$ deux poids raisonnables respectivement sur
$B_+(E_1)$ et
$B_+(E_2)$. On peut poser pour $u = \sum_1^n x_i \otimes y_i \in E_1
\otimes
E_2$
$$\gamma(u) = \inf \{ \gamma_1(\sum x_i \otimes \adh {x_i})^{1/2}
\gamma_2
(\sum y_i \otimes \adh {y_i})^{1/2}\}\leqno (1)$$
o\`u l'infimum porte sur toutes les repr\'esentations possibles de $u$.
On v\'erifie
alors ais\'ement que $\gamma$ est une norme sur $E_1 \otimes E_2$. On
notera $E_1
\hat \otimes_\gamma E_2$ l'espace compl\'et\'e associ\'e. On peut
montrer que la
norme duale de $\gamma$ est essentiellement du m\^eme type. Plus
pr\'ecis\'ement, on a

\proclaim THEOREME 4. Soient $E_1,E_2$ deux espaces de dimension finies
munis
de poids raisonnables $\gamma_1$ sur $B_+(E_1)$ et $\gamma_2$ sur
$B_+(E_2)$.
On notera $\gamma_1$ (resp. $\gamma_2)$ la norme \'etendant $\gamma_1$
sur $E_1
\otimes \adh {E_1}$ (resp. $E_2 \otimes \adh {E_2})$ et $\gamma_1^*$
(resp.
$\gamma_2^*)$ la norme duale sur $E_1^* \otimes \adh {E_1^*}$ (resp.
$E_2^*
\otimes \adh {E_2^*})$. Soit $\gamma$ la norme associ\'ee sur $E_1
\otimes E_2$
comme en (1) ci-dessus. Alors sa norme duale $\gamma^*$ coincide avec
la norme
associ\'ee comme en (1) ci-dessus avec les normes $\gamma_1^*$ et
$\gamma_2^*$,
c'est-\`a-dire que pour tout $v$ dans $E_1^* \otimes E_2^*$ on a
$$\gamma^*(v) = \inf \{\gamma_1^* (\sum \xi_i \otimes \adh
{\xi_i})^{1/2}
\gamma_2^* (\sum \eta_1 \otimes \adh {\eta_i})^{1/2}\}$$
o\`u l'infimum porte sur toutes les repr\'esentations possibles de la
forme $v =
\sum_1^n \xi_i \otimes \eta_i$ avec $\xi_i \in E_1^*,\  \eta_i \in
E_2^*$.

On peut voir le th\'eor\`eme 2 comme un cas particulier du th\'eor\`eme
4, en prenant
$\gamma_1$ et $\gamma_2$ \'egales \`a la norme du produit tensoriel
minimal
(=spatial). On notera que l'on peut d\'eduire du th\'eor\`eme 2 (ou du
th\'eor\`eme 4)
une description de la norme tensorielle $\gamma_{oh}^*$
comme suit.
Soient $E_1,E_2$ deux espaces d'op\'erateurs. Pour tout $u\in
E_1\otimes E_2$
 posons
$$\gamma_{oh}^*(u)=\sup\{ |<u,v>|\quad  |\quad v\in E_1^*\otimes
E_2^*,\
\gamma_{oh}(v)\leq 1\}.$$
Soit $E$ un espace d'op\'erateurs arbitraire et soient $u,v\in
(E\om\adh E)^*$
on note $u\le v$ si $u$ et $v$ sont sym\'etriques et si $u(x,x)\le
v(x,x)$ pour tout $x$ dans $ E$. Comme
d'habitude un \'el\'ement de $E\otimes\adh E$ peut aussi \^etre
consid\'er\'e comme un
\'el\'ement de $(E^*\om\adh{ E}^*)^*$.
On peut alors poser
$$\forall u\in E\om\adh E\quad {\rm avec}\quad u\ge 0
\quad \alpha(u)=\inf\{I(v) |\quad v\in (E^*\om\adh{ E}^*)^*,\quad u\le
v\}.$$
On a alors
\proclaim COROLLAIRE 5. Soient $E_1,E_2$ deux espaces d'op\'erateurs.
 Pour tout $u\in E_1\otimes E_2$ on a
$$\gamma_{oh}^*(u)= \inf \{\alpha(\sum x_i
\otimes \adh {x_i})^{1/2} \alpha(\sum y_i \otimes \adh {y_i})^{1/2}\}$$
o\`u l'infimum porte sur toutes les repr\'esentations $u = \sum_1^n x_i
\otimes
  {y_i},\quad
x_i \in E_1,\quad y_i \in E_2.$

Parmi les cons\'equences notons que pour tout $u\in E_1\otimes E_2$ on
a
$\gamma_{oh}^*(u)\geq \gamma_{oh}(u)$. De plus, si $E_2$ est suppos\'e
r\'eflexif
(pour simplifier l'\'enonc\'e) on a
$$\gamma_{oh}^*(u)=\inf\{ \pi_{2,oh}(B) \pi_{2,oh}(A^*) \}$$
o\`u l'infimum porte sur toutes les repr\'esentations  possibles de $u$
de la forme
$u=AB$ avec $B: E_1^*\to \ell_2$  et $A:\ell_2\to E_2$ tels que
  $B$ et $A^*$ sont $(2,oh)-$sommants.
Enfin pour tout op\'erateur $v:E_1 \to OH$ on a
$\pi_{2,oh}(v)\geq \|v\|_{cb}.$

 Signalons
qu'on peut montrer en suivant des id\'ees bien connues que, si la
norme $\gamma$ est comme en (1), l'application canonique de $E_1 \hat
\otimes_\gamma E_2$ dans le produit tensoriel injectif $E_1 \buildrel
{\vee}\over {\otimes} E_2$ est toujours injective quels que soient les
espaces
de Banach $E_1$ et $E_2$.

Les id\'ees pr\'ec\'edentes permettent de d\'evelopper des notions
d'espace
d'op\'erateurs de type 2 ou de cotype 2. Soit $E$ un espace
d'op\'erateurs, $F$ un
espace de Banach et $u : F \to E$ un op\'erateur. On notera
$\Pi_{2,oh}^* (F,E)$
la classe des op\'erateurs $w : F \to E$ admettant une factorisation de
la forme
$F \to OH \to E$ avec $A : OH \to E$ compl\`etement born\'e et $B : F
\to OH$ tel
que pour une base orthonorm\'ee $(T_n)$ de $OH$ on a $\sum \N
{B^*(T_n)}^2 <
\infty$. On posera $\pi_{2,oh}^*(w) = \inf \{\N {A}_{cb} (\sum \N
{B^*(T_n)}^2)^{1/2}\}$ o\`u l'infimum porte sur toutes les
factorisations
possibles de $w$.

Soit $E$ un espace d'op\'erateur et $F$ un espace de Banach. On posera
pour $v
\in E \otimes F$
$$d_{2,oh}(v) = \inf \{\N {\sum x_i \otimes \adh {x_i}}^{1/2} (\sum \N
{y_i}^2)^{1/2}\}$$
o\`u l'infimum porte sur toutes les repr\'esentations
$$v = \sum_1^n x_i \otimes y_i,\quad x_i \in E,\quad y_i \in F.$$
Soit $E \hat \otimes_{d_{2,oh}} F$ le compl\'et\'e associ\'e. Il est
facile
de v\'erifier (en suivant des id\'ees bien connues des sp\'ecialistes
des espaces de
Banach) que $(E \hat \otimes_{d_{2,oh}} F)^*$ s'identifie
isom\'etriquement \`a
l'espace $\Pi_{2,oh} (E,F^*)$. De mani\`ere \'equivalente si l'on
associe \`a $v \in E
\otimes F$ un op\'erateur $\tilde v : F^* \to E$, on voit ais\'ement
que
$$d_{2,oh} (v) = \inf \{ \N {A}_{cb} (\sum \N {B^* T_n}^2)^{1/2}\}$$
o\`u l'infimum porte sur toutes les factorisations $\tilde v = AB$ avec
$B : F^*
\to OH, A : OH \to E$ et o\`u $(T_n)$ est une base orthonormale fix\'ee
de $OH$.
Soit alors $E \subset B(H)$ un espace d'op\'erateurs. Soit $(e_i)$ la
base
canonique de $\ell_2$.
Soit $x = (x_i)$ une suite finie dans $E$, on note $u_x = \sum x_i
\otimes
e_i$ l'op\'erateur de $E^*$ dans $\ell_2$ d\'efini par $u_x(\xi) = \sum
\xi(x_i)
e_i$. Soit $(g_i)$ une suite de variables al\'eatoires gaussiennes
ind\'ependantes
standard sur un espace de probabilit\'e $(\Omega,{\cal A},P)$. On
notera $L^2(E)$
l'espace $L^2(\Omega,{\cal A},P ; E)$.

Nous dirons qu'un espace d'op\'erateurs $E$ est de $OH$-type 2 (resp.
$OH$-cotype
2) s'il existe une constante $C$ telle que pour toute suite finie
$(x_i)$ dans
$E$ on a
$$\N {\sum g_i x_i}_{L^2(E)} \leq C \pi_{2,oh} (\sum x_i \otimes
e_i).$$
$${\rm (resp.} d_{2,oh} (\sum x_i \otimes e_i) \leq C \N {\sum g_i
x_i}_{L^2(E)}).$$ On peut v\'erifier par exemple (voir [6] pour des
d\'efinitions
pr\'ecises)  que l'espace $R\cap C$  est de $OH$-type 2 et que $R+C$
est
de $OH$-cotype 2. On obtient alors ais\'ement l'analogue pour les
espaces
d'op\'erateurs d'une s\'erie de th\'eor\`emes classiques en ``th\'eorie
locale" des
espaces de Banach. Par exemple, on peut \'etendre un th\'eor\`eme de
Kwapie\'n (voir
[7]) : tout op\'erateur  born\'e d'un espace $E$ de $OH$-type 2 dans un
espace $F$
de $OH$-cotype 2 se factorise par $OH$ (et a fortiori est
compl\`etement born\'e).

Notons qu'un espace de $OH$-type 2 (resp. $OH$-cotype 2) est a fortiori
de type
2 (cotype 2) au sens usuel puisque pour tout \'el\'ement $\sum_1^n x_i
\otimes y_i$
de $E \otimes F$ les op\'erateurs associ\'es $u : E^* \to F$ et $u^* :
F^* \to E$
v\'erifient $\pi_{2,oh}(u) \leq \pi_2(u)$ et $\pi_2(u^*) \leq d_{2,oh}
(\sum x_i
\otimes y_i)$.

Je remercie David Blecher et Vern Paulsen pour des conversations
stimulantes
sur le sujet  de cette note et de la note pr\'ec\'edente
[6].\vfill\eject
{\bf References}
\item{[1]} D. Blecher and V. Paulsen. Tensor products of operator
spaces.
 J. Funct. Anal. 99 (1991) 262-292.

 \item{[2]} E. Effros and Z.J. Ruan. A new approach to operator
 spaces.
 Canadian Math. Bull.34(1991) 329-337.

  \item{[3]} E. Effros and Z.J. Ruan. Self-duality for the Haagerup
  tensor
product and Hilbert space factorization, J. Funct. Anal. (A para\^\i
tre).

\item{[4]} E. Effros and Z.J. Ruan. Mapping spaces and liftings for
operator
spaces. A para\^\i tre.

\item {[5]} A. Grothendieck  R\'esum\'e de la th\'eorie m\'etrique des
produits tensoriels topologiques.  Boll.. Soc. Mat.
S$\tilde{a}$o-Paulo  8
(1956), 1-79.

\item {[6]} G.Pisier. Espace de Hilbert d'op\'erateurs et interpolation
complexe.
C. R. Acad. Sci. Paris

\item {[7]} $\underline{\hskip1.5in}$. Factorization of linear
operators and the Geometry of Banach spaces.  CBMS
(Regional conferences of the A.M.S.)    60, (1986),
Reprinted with corrections 1987.

\item{[8]} G. Pisier. Factorization of operator valued analytic
functions.
 Advances in Math. 93 (1992) 61-125.

 \medskip

Equipe d'Analyse

Universit\'e Paris 6

Bo\^\i te 186, Place Jussieu

75252 Paris Cedex 05, FRANCE

et

Mathematics Department

Texas A. and M. University

College Station, Texas 77843, USA

 \end

consid\'er\'e comme un op\'erateur de $E_1^*$ dans $E_2$
\proclaim COROLLAIRE 5 Soient $E_1,E_2$ deux espaces d'op\'erateurs.
 Pour tout $u\in E_1\otimes E_2$ on a
$$\gamma_{oh}^*(u)=\inf\{ \pi_{2,oh}(B) \pi_{2,oh}(A^*) \}$$
o\`u l'infimum porte sur toutes les repr\'esentations  possibles de $u$
de la forme
$u=AB$ avec $B: E_1^*\to \ell_2$  et $A:\ell_2\to E_2$ tels que
  $B$ et $A^*$ sont $(2,oh)-$sommants.

\end
\magnification\magstep1
 \baselineskip = 18pt
 
 \overfullrule = 0pt
 
\pageno=44

 \centerline{\bf References}

 \item{[BL]} J. Bergh and J. L\"ofstr\"om. Interpolation spaces. An
 introduction. Springer Verlag, New York. 1976.

 \item{[BP]} D. Blecher and V. Paulsen. Tensor products of operator
 spaces.
 J. Funct. Anal. (1990).

 \item{[BS]} D. Blecher and R. Smith. The dual of the Haagerup tensor
 product. (Preprint) 1990.

 \item{[B1]} D. Blecher. Tensor products of operator spaces II.
 (Preprint)
 1990. To appear in Pacific J. Math.

 \item{[B2]} D. Blecher. The standard dual of an operator space.

 \item{[Ca]} A. Calder\'on. Intermediate spaces and interpolation, the
 complex method. Studia Math. 24 (1964) 113-190.

\item{ [DCH]} J. de Canni\`ere   and  U. Haagerup.
Multipliers of the Fourier algebras of some simple Lie
groups and their discrete subgroups.  Amer. J. Math.   107
(1985), 455-500.

 \item{[EE]} E. Effros and R. Exel. On multilinear double commutant
 theorems.  Operator algebras and applications, Vol.~1, London Math.
 Soc.
 Lecture Notes Series 135 (19\ \ ) 81-94.

 \item{[EK]} E. Effros and A. Kishimoto. Module maps and
 Hochschild-Johnson cohomology. Indiana Univ. Math. J. 36 (1987)
 257-276.

 \item{[ER1]}  E. Effros and Z.J. Ruan. A new approach to operator
 spaces.
 Canadian Math. Bull. 34(1991) 329-337.

 \item{[ER2]} $\underline{\hskip1.5in}$. A new approach to operators
 spaces.

 \item{[ER3]} $\underline{\hskip1.5in}$. On the abstract
 characterization of
 operator spaces.

 \item{[ER4]} $\underline{\hskip1.5in}$. Self duality for the Haagerup
 tensor product and Hilbert space factorization. (Preprint) 1990.

 \item{[ER5]} $\underline{\hskip1.5in}$. Recent development in operator
 spaces.

\item{[ER6]} E. Effros and Z.J. Ruan. On matricially
normed spaces. Pacific
 J. Math. 132 (1988) 243-264.

\item {[GK]} J. Grosberg and M. Krein.  Sur la
d\'ecomposition des fonctionnelles en composantes
positives.  Doklady Akad. Nauk SSSR {\bf 25} (1939),
723-726.

\item{[G]} A. Grothendieck.  R\'esum\'e de la th\'eorie
m\'etrique des produits tensoriels topologiques.  Boll..
Soc. Mat. S$\tilde{a}$o-Paulo  8 (1956), 1-79.

 \item{[HP1]} U. Haagerup and G. Pisier. Factorization of analytic
 functions
 with values in non-commutative $L^1$-spaces and applications. Canadian
 J.
 Math. 41 (1989) 882-906.

 \item{[HP2]} U. Haagerup and G. Pisier. Bounded linear
operators between
 $C^*$-algebras. To appear.

 \item{[H1]} U. Haagerup. Injectivity and decomposition of completely
 bounded maps in ``Operator algebras and their connection with Topology
 and
 Ergodic Theory''. Springer Lecture Notes in Math. 1132 (1985) 91-116.

 \item{[H2]} $\underline{\hskip1.5in}$. An example of a non-nuclear
 $C^*$-algebra which
 has the metric approximation property. Inventiones Math. 50 (1979)
 279-293.

 \item{[H3]} $\underline{\hskip1.5in}$. Decomposition of completely
 bounded
 maps on operator algebras. Unpublished manuscript. Sept. 1980.

\item {[ KR ]} R. Kadison  and J. Ringrose.    Fundamentals
of the theory of operator algebras, Vol. II, Advanced
Theory,    Academic Press, New-York      1986.

 \item{[K]} O. Kouba. Interpolation of injective or projective tensor
 products of Banach spaces. J. Funct. Anal. 96 (1991) 38-61.

 \item{[L]} C. Le Merdy. Analytic factorizations and completely bounded
 maps. Preprint. June 1992.

 \item{[Pa1]} V. Paulsen. Representation of Function algebras, Abstract
 operator spaces and Banach space Geometry. Preprint.

 \item{[PS]} V. Paulsen and R.Smith. Multilinear maps and
tensor norms on operator systems. J. Funct. Anal. 73
(1987) 258-276.

 \item{[P1]} G. Pisier. Factorization of operator valued analytic
 functions.
 Advances in Math. 93 (1992) 61-125.

 \item{[P2]} $\underline{\hskip1.5in}$. Factorizations of operator
 valued
 analytic functions and complex interpolation. Feitschrift in honor of
 I.~Piatetski-Shapiro (edited by Gelbard, Howe and Sarnek) Part II,
 pp.~197-220, Weizmann Science Press of Israel, IMCP, 1990.

\item{[P3]} $\underline{\hskip1.5in}$. Remarks on
complemented subspaces of $C^*$-algebras. Proc. Roy. Soc.
Edinburgh, 121A (1992) 1-4.

\item {[P4]} $\underline{\hskip1.5in}$. Factorization of linear
operators and the Geometry of Banach spaces.  CBMS
(Regional conferences of the A.M.S.)    60, (1986),
Reprinted with corrections 1987.

\item{[RW]} A.R. Robertson  and S.Wasserman.

\item{[Ru]} Z.J. Ruan. Subspaces of $C^*$-algebras. J. Funct. Anal. 76
 (1988) 217-230.

\item{[T]} J.Tomiyama. On the projection of norm one in
$W^*$-algebras. Proc. Japan Acad. 33 (1957) 608-612.
  \medskip

Equipe d'Analyse

Universit\'e Paris 6

Bo\^\i te 186, Place Jussieu

75252 Paris Cedex 05, FRANCE

et

Mathematics Department

Texas A. and M. University

College Station, Texas 77843, USA

 \end

\end